\documentclass{amsart}
\usepackage{latexsym}
\usepackage{amsfonts}
\usepackage{amssymb}
\usepackage{graphicx}

\newtheorem{theorem}{Theorem}[section]
\newtheorem{corollary}[theorem]{Corollary}
\newtheorem{lemma}[theorem]{Lemma}
\newtheorem{proposition}[theorem]{Proposition}
\theoremstyle{remark}
\newtheorem{remark}{Remark}[section]
\numberwithin{equation}{section}

\newcommand{\g}{\mathfrak{g}}
\newcommand{\ak}{\mathfrak{k}}
\newcommand{\am}{\mathfrak{m}}
\newcommand{\RR}{\mathbb{R}}

\begin{document}

\title[Curve Selection Lemma and conjugacy classes in Lie groups]
{Curve Selection Lemma for semianalytic sets and conjugacy classes
of finite order in Lie groups}

\author{Jinpeng An}
\address{School of mathematical science, Peking University,
 Beijing, 100871, P. R. China }
\email{anjinpeng@pku.edu.cn}

\author{Zhengdong Wang}
\address{School of mathematical science, Peking University,
 Beijing, 100871, P. R. China}
\email{zdwang@pku.edu.cn}

\thanks{This work is supported by the 973 Project
Foundation of China ($\sharp$TG1999075102).}

\keywords{semianalytic set, Curve Selection Lemma, Lie group,
conjugacy class.}

\subjclass[2000]{Primary 14P15; Secondary 22E15.}

\begin{abstract}
Using a strong version of the Curve Selection Lemma for real
semianalytic sets, we prove that for an arbitrary connected Lie
group $G$, each connected component of the set $E_n(G)$ of all
elements of order $n$ in $G$ is a conjugacy class in $G$. In
particular, all conjugacy classes of finite order in $G$ are
closed. Some properties of connected components of $E_n(G)$ are
also given.
\end{abstract}

\maketitle


\section{Introduction}

The interest of real semianalytic geometry (or, more generally,
real subanalytic geometry) lies not only on its own right but also
on its applications in many fields, such as calculus of variations
and control theory. In this paper we give an application of real
semianalytic geometry to the problem of the closedness of
conjugacy classes of finite order in Lie groups. The availability
of real semianalytic geometry owes to the fact that each Lie group
admits a unique real analytic structure such that multiplication,
inversion, the exponential map, and homomorphisms are real
analytic (see, for example, Varadarajan \cite{Va}, Section 2.11).
Using a strong version of the Curve Selection Lemma for real
semianalytic sets, we will prove the following conclusion.

\begin{theorem}\label{T:main}
Let $G$ be a connected Lie group, $n$ a positive integer. Denote
$E_n(G)=\{g\in G|g^n=e\}$. Then each connected component of
$E_n(G)$ is a conjugacy class of $G$.
\end{theorem}

Since $E_n(G)$ is a closed subset of $G$, its connected components
are closed in $G$. So Theorem \ref{T:main} implies immediately
that all conjugacy classes of order $n$ in $G$ are closed.

It is well known that semisimple conjugacy classes in algebraic
groups are Zariski closed (see, for example, Humphreys \cite{Hu}),
and a conjugacy class of finite order in an algebraic group is
necessary semisimple. So Theorem \ref{T:main} may be viewed as a
generalization of this result to Lie groups.

In Section 2, we will present the version of the Curve Selection
Lemma that will be needed. Based on it, the proof of Theorem
\ref{T:main} will be given in Section 3.

As corollaries of Theorem \ref{T:main}, some properties of the
connected components of the set $E_n(G)$ will be given in Section
4. Denote the set of all connected components of $E_n(G)$ by
$\mathcal{E}_n(G)$. The most fundamental property of
$\mathcal{E}_n(G)$ is that each element of $\mathcal{E}_n(G)$ is a
closed submanifold of $G$, which is an obvious corollary of
Theorem \ref{T:main}. Of course, different elements of
$\mathcal{E}_n(G)$ may have different dimensions which, as the
dimensions of conjugacy classes, can be easily computed. Another
corollary is that
$\mathcal{E}_n(G)\cap\mathcal{E}_m(G)=\mathcal{E}_{(n,m)}(G)$ for
any positive integers $n,m$, where $(n,m)$ is the greatest common
divisor of $n$ and $m$. The following property will also be
proved. Let $K$ be a maximal compact subgroup of $G$, and let
$i:K\rightarrow G$ be the inclusion map. We will show that the
natural induced map $i_*:\mathcal{E}_n(K)\rightarrow
\mathcal{E}_n(G)$ is a bijection. In particular,
$\mathcal{E}_n(G)$ is a finite set, and the intersection of each
connected component of $E_n(G)$ with $K$ is a connected component
of $E_n(K)$.


\section{A version of the Curve Selection Lemma for semianalytic sets}

In this section we present a strong version of the Curve Selection
Lemma for real semianalytic sets, which will be needed in the
proof of Theorem \ref{T:main}. For a real semianalytic set $X$,
the sets of regular points and singular points in $X$ are denoted
by $r(X)$ and $s(X)$, respectively. The following Wing Lemma is a
key ingredient.

\begin{proposition}\label{P:wl}
\textbf{(Wing Lemma \cite{Lo1,Wa})} Let $m<n$ be non-negative
integers, and identify $\RR^m$ with
$\{(x_1,\cdots,x_n)\in\RR^n|x_{m+1}=\cdots=x_n=0\}$. Let
$U\subset\RR^n$ be an open neighborhood of $0\in\RR^n$, and let
$X\subset U\cap(\RR^{n}-\RR^m)$ be a semianalytic set. Suppose
$U\cap\RR^m\subset\overline{X}$. Then there exist semianalytic
sets $W,Z$ of $\RR^n$ such that\\
(i) $W\subset X$, $U\cap\RR^m\subset\overline{W}$, $\dim(W)=m+1$;\\
(ii) $Z\subset U\cap\RR^m$, $Z$ is closed and nowhere dense in
$U\cap\RR^m$, $\dim(Z)<m$;\\
(iii) For each $p\in (U\cap\RR^m)-Z$, there is an open
neighbourhood $V\subset U$ of $p$ in $\RR^n$ such that for each
connected component $W_r$ of $W\cap V$, $W_r\cup(V\cap\RR^m)$ is
an $(m+1)$-dimensional $C^1$ submanifold of $\RR^n$ with boundary
$V\cap\RR^m$, and there exists a $C^1$ function
$y^{r}:W_r\cup(V\cap\RR^m)\rightarrow \RR$ such that
$(x_1,\cdots,x_m,y^{r})$ forms a $C^1$ coordinate system of
$W_r\cup(V\cap\RR^m)$.
\end{proposition}

\begin{proposition}\label{P:curve}
\textbf{(Curve Selection Lemma)} Let $X$ be a semianalytic set in
a real analytic manifold $M$. Suppose $s(X)\neq\emptyset$. Then
there exists a semianalytic set $Z\subset s(X)$ which is closed
and nowhere dense in $s(X)$ with $\dim(Z)<\dim(s(X))$ such that
for each $p\in s(X)-Z$, there exists a $C^1$ curve
$\gamma:[0,1)\rightarrow M$ such that $\gamma(0)=p$,
$\gamma((0,1))\subset r(X)$, and $\gamma'(0)$ does not belong to
the tangent space $T_pN$ of any analytic immersed submanifold $N$
of $M$ with $p\in N\subset s(X)$.
\end{proposition}

\begin{proof}
Choose a locally finite family of analytic coordinate charts
$\{(U^i,x^i_1,\cdots,x^i_n)|i\in I\}$ such that
$r(s(X))\subset\bigcup_{i\in I}U^i$ and $U^i\cap s(X)=\{p\in
U^i|x^i_{m_i+1}(p)=\cdots x^i_n(p)=0\}$, where $n=\dim(M),
m_i\leq\dim(s(X))$ (such a family of charts always exists since
$r(s(X))$ is open in $s(X)$). By the Wing Lemma, for each $i\in
I$, there exists a semianalytic set $Z^i\subset U^i\cap s(X)$
which is closed and nowhere dense in $U^i\cap s(X)$ with
$\dim(Z^i)<m_i$ such that for each $p\in (U^i\cap s(X))-Z^i$,
there exists a $C^1$ curve $\gamma:[0,1)\rightarrow U^i$ such that
$\gamma(0)=p$, $\gamma((0,1))\subset r(X)\cap U^i$, and
$\gamma'(0)\notin T_p(s(X))$ (note that $T_p(s(X))$ is
well-defined since $p\in r(s(X))$). In fact, the $y^{r}$-axis in
Proposition \ref{P:wl} may serve as the curve $\gamma$. Now the
semianalytic set $Z=s(s(X))\cup\bigcup_{i\in I}Z^i$ satisfies our
requirement, since $s(s(X))$ is nowhere dense in $s(X)$.
\end{proof}

\begin{remark}
The elimination of a semianalytic subset $Z$ of small dimension is
necessary in Proposition \ref{P:curve}. For example, consider the
analytic subset $X$ of $\RR^3$ defined by $(y^2+z^2)^2-4x^4z^2=0$.
Then the intersection of $X$ with the plane $x=a$ is two tangent
circles with radii $a^2$ and tangent point $(a,0,0)$. The set of
singular points $s(X)$ is the whole $x$-axis. But for each $C^1$
curve starting with $(0,0,0)$ and being contained in $X$, its
tangent vector at $(0,0,0)$ belongs to the $x$-axis.
\end{remark}


\section{Conjugacy classes of finite order in Lie groups}

With the aid of the Curve Selection Lemma proved in the previous
section, we can prove the main theorem of this paper now. For a
Lie group $G$, denote $E_n(G)=\{g\in G|g^n=e\}$. First we give
several lemmas.

\begin{lemma}\label{L:dimOg}
Let $G$ be a Lie group, $g\in G$. Let $C_g$ be the conjugacy class
containing $g$. Then as an immersed submanifold of $G$, its
tangent space at $g$ is
$$
T_gC_g=(dr_g)_e(\mathrm{Im}(1-Ad(g))),
$$
where $r_g$ is the right translation on $G$ by $g$.
\end{lemma}

\begin{proof}
A direct computation.
\end{proof}

\begin{lemma}\label{L:dimEn}
Let $G$ be a Lie group, $g\in E_n(G)$. Let
$\gamma:[0,1)\rightarrow G$ be a $C^1$ curve such that
$\gamma(0)=g$, $\gamma((0,1))\subset E_n(G)$. Then
$$
\gamma'(0)\in(dr_g)_e(\ker(1+Ad(g)+\cdots+Ad(g^{n-1}))).
$$
\end{lemma}

\begin{proof}
Let $\alpha(t)=\gamma(t)g^{-1}$. Then
$$
e=\gamma(t)^n=(\alpha(t)g)^n=
\alpha(t)(g\alpha(t)g^{-1})\cdots(g^{n-1}\alpha(t)g^{-(n-1)}).
$$
Applying $\frac{d}{dt}\big|_{t=0}$ two both sides of the above
identity, we get
$$
(1+Ad(g)+\cdots+Ad(g^{n-1}))(\alpha'(0))=0.
$$
Note also that $\gamma'(0)=(dr_g)_e(\alpha'(0))$, the lemma is
proved.
\end{proof}

\begin{lemma}\label{L:dim}
Let $G$ be a Lie group, $g\in E_n(G)$. Then
$$
\ker(1+Ad(g)+\cdots+Ad(g^{n-1}))=\mathrm{Im}(1-Ad(g))
$$
on the Lie algebra $\g$ of $G$.
\end{lemma}

\begin{proof}
First we have $(1+Ad(g)+\cdots+Ad(g^{n-1}))\circ(1-Ad(g))=0$. So
$$
\mathrm{Im}(1-Ad(g))\subset\ker(1+Ad(g)+\cdots+Ad(g^{n-1})).
$$
Hence to prove the lemma, it is sufficient to show that
$$
\mathrm{rank}(1-Ad(g))+\mathrm{rank}(1+Ad(g)+\cdots+Ad(g^{n-1}))\geq\dim(\g).
$$
We claim that
$$
\ker(1-Ad(g))\cap\ker(1+Ad(g)+\cdots+Ad(g^{n-1}))=\{0\}.
$$
If this is true, then the inequality about ranks follows
immediately. Now we verify the claim. In fact, if $X\in\g$ belongs
to the left hand side of the above equality, then $(1-Ad(g))X=0$,
which means that $Ad(g)X=X$. So
$0=(1+Ad(g)+\cdots+Ad(g^{n-1}))X=nX$. This shows $X=0$. The claim
is proved.
\end{proof}

\begin{theorem}\label{T:En}
Let $G$ be a connected Lie group, $n$ a positive integer. Then
each connected component of $E_n(G)$ is a conjugacy class of $G$.
\end{theorem}

\begin{proof}
Endow $G$ with the unique real analytic structure such that group
operations are real analytic, then $E_n(G)$ is an analytic subset
of $G$. Let $E^i$ be a connected component of $E_n(G)$, then $E^i$
is semianalytic. Since $E_n(G)$ is invariant under the adjoint
action of $G$ and $G$ is connected, $E^i$ is invariant under the
adjoint action of $G$. Hence $r(E^i)$ and $s(E^i)$ are also
invariant under the adjoint action of $G$. We claim that
$s(E^i)=\emptyset$. For otherwise, by Proposition \ref{P:curve},
there exists a $g\in s(E^i)$ and a $C^1$ curve
$\gamma:[0,1)\rightarrow G$ such that $\gamma(0)=g$,
$\gamma((0,1))\subset r(E^i)$, and $\gamma'(0)\notin T_gN$ for any
analytic immersed submanifold $N$ of $G$ with $g\in N\subset
s(E^i)$. In particular, $\gamma'(0)\notin T_gC_g$, where $C_g$ is
the conjugacy class containing $g$, which is an analytic immersed
submanifold of $G$. But by Lemma \ref{L:dimOg}, Lemma
\ref{L:dimEn}, and Lemma \ref{L:dim}, we have $\gamma'(0)\in
T_gC_g$, a contradiction. This shows $s(E^i)=\emptyset$, that is,
$E^i$ is an analytic submanifold of $G$.

The proof of Lemma \ref{L:dimEn} in fact shows
$$
T_gE^i=(dr_g)_e(\ker(1+Ad(g)+\cdots+Ad(g^{n-1})))
$$
for each $g\in E^i$. By Lemma \ref{L:dimOg} and Lemma \ref{L:dim},
$T_gE^i=T_gC_g$. This means that the infinitesimal transformations
of $\g$ on $E^i$ induced by the adjoint action generate $T_gE^i$
at each $g\in E^i$. So the adjoint action of $G$ acts transitively
on $E^i$. In other word, we have $E^i=C_g$. The theorem is proved.
\end{proof}

Let $n$ be a positive integer. An element $g$ in a Lie group is of
\emph{order $n$} if $g^n=e$. A conjugacy class in a Lie group is
of \emph{order $n$} if it contains an element of order $n$, in
which case all elements in this conjugacy class are of order $n$.
The following corollary is obvious.

\begin{corollary}
Let $G$ be a connected Lie group, then all conjugacy classes of
order $n$ in $G$ are closed for each positive integer $n$.\qed
\end{corollary}

\begin{remark}
By an argument similar to the proof of Theorem \ref{T:En}, one can
easily show that for a Lie group $G$ with finite many components,
each conjugacy class in $G$ of order $n$ is a union of finite many
connected components of $E_n(G)$, which is closed. We leave the
details to the reader.
\end{remark}

The following proposition gives a description of the subsets
$E_n(G)$.

\begin{proposition}
Let $G$ be a connected Lie group. Then
$\overline{\bigcup_{n=1}^\infty
E_n(G)}=\overline{\bigcup_{K\in\mathcal{K}(G)}K}$, where
$\mathcal{K}(G)$ denotes the set of all maximal compact subgroups
of $G$. In particular, if $G$ is compact, then
$\bigcup_{n=1}^\infty E_n(G)$ is dense in $G$.
\end{proposition}

\begin{proof}
Let $n$ be a positive integer and $g\in E_n$. Then the set
$\{e,g,\cdots,g^{n-1}\}$ is a finite (hence compact) subgroup of
$G$. By Theorem 3.1 of Chapter XV in Hochschild \cite{Ho}, there
exists a maximal compact subgroup $K$ of $G$ containing $g$. So
$E_n\subset\bigcup_{K\in\mathcal{K}(G)}K$, and hence
$\overline{\bigcup_{n=1}^\infty
E_n(G)}\subset\overline{\bigcup_{K\in\mathcal{K}(G)}K}$.

Conversely, let $K$ be a maximal compact subgroup $K$ of $G$,
which is necessary connected (also by the loc. cit. Theorem in
\cite{Ho}). Let $T$ be a maximal torus of $K$. Then it is obvious
that $T=\overline{\bigcup_{n=1}^\infty
E_n(T)}\subset\overline{\bigcup_{n=1}^\infty E_n(G)}$. But $K$ is
the union of all maximal tori of $K$. So
$K\subset\overline{\bigcup_{n=1}^\infty E_n(G)}$, and hence
$\overline{\bigcup_{K\in\mathcal{K}(G)}K}\subset\overline{\bigcup_{n=1}^\infty
E_n(G)}$. This proves the proposition.
\end{proof}


\section{Some properties of connected components of $E_n(G)$}

In this section we present some properties of connected components
of the subset $E_n(G)$. To simply notations, We denote the set of
all connected components of $E_n(G)$ by $\mathcal{E}_n(G)$.

\begin{corollary}
Let $G$ be a connected Lie group, $n$ a positive integer. Then
each element of $\mathcal{E}_n(G)$ is a closed submanifold of $G$.
\end{corollary}

\begin{proof}
Let $E^i$ be an element of $\mathcal{E}_n(G)$, that is, a
connected component of $E_n(G)$. By Theorem \ref{T:En}, $E^i$ is a
conjugacy class in $G$, so it is an immersed submanifold of $G$.
But as a connected component of $E_n(G)$, $E^i$ is closed in
$E_n(G)$, hence closed in $G$. This proves $E^i$ is a closed
submanifold of $G$.
\end{proof}

It is obvious that $E_n(G)\cap E_m(G)=E_{(n,m)}(G)$, where $(n,m)$
is the greatest common divisor of $n$ and $m$. But without the
help of Theorem \ref{T:En}, the following conclusion is not so
obvious.

\begin{corollary}\label{C:nm}
Let $G$ be a connected Lie group, and let $n,m$ be positive
integers. Then
$\mathcal{E}_n(G)\cap\mathcal{E}_m(G)=\mathcal{E}_{(n,m)}(G)$. In
particular, if $n$ divides $m$, then
$\mathcal{E}_n(G)\subset\mathcal{E}_m(G)$.
\end{corollary}

\begin{proof}
Suppose $E^i\in\mathcal{E}_{(n,m)}(G)$. By Theorem \ref{T:En},
$E^i$ is a conjugacy class. Note that $E^i\subset E_n(G)$, also by
Theorem \ref{T:En}, $E^i\in\mathcal{E}_n(G)$. Similarly,
$E^i\in\mathcal{E}_m(G)$. So we have
$\mathcal{E}_{(n,m)}(G)\subset\mathcal{E}_n(G)\cap\mathcal{E}_m(G)$.

Conversely, suppose $E^j\in\mathcal{E}_n(G)\cap\mathcal{E}_m(G)$.
Then $E^j$ is a conjugacy class contained in $E_{(n,m)}(G)$, and
then is an element of $\mathcal{E}_{(n,m)}(G)$. So we also have
$\mathcal{E}_n(G)\cap\mathcal{E}_m(G)\subset\mathcal{E}_{(n,m)}(G)$.
This proves the corollary.
\end{proof}

\begin{remark}
Since $E_n(G)\cap E_m(G)=E_{(n,m)}(G)$, Corollary \ref{C:nm} in
fact says $\mathcal{E}_n(G)\cap\mathcal{E}_m(G)$ coincides with
the set of all connected components of $E_n(G)\cap E_m(G)$.
\end{remark}

Let $H$ be another Lie group and let $f:H\rightarrow G$ be a
homomorphism of Lie groups. Then $f$ induces a map
$f_*:\mathcal{E}_n(H)\rightarrow \mathcal{E}_n(G)$ in the natural
way. In the case that $K$ is a maximal compact subgroup of $G$ and
$i:K\rightarrow G$ is the inclusion map, the induced map
$i_*:\mathcal{E}_n(K)\rightarrow \mathcal{E}_n(G)$ is a bijection.
To prove this, we need the following lemma.

\begin{lemma}\label{L:max}
Let $G$ be a connected Lie group, $K$ a maximal compact subgroup
of $G$. Then two elements of $K$ are conjugate in $K$ if and only
if they are conjugate in $G$.
\end{lemma}

\begin{proof}
The ``only if" part is obvious. We prove the ``if" part.

By Theorem 3.1 of Chapter XV in \cite{Ho}, there exist some linear
subspaces $\am_1,\cdots,\am_r$ of the Lie algebra $\g$ of $G$ such
that\\
(1) $\g=\ak\oplus\am_1\oplus\cdots\oplus\am_r$, where $\ak$ is the
Lie algebra of $K$;\\
(2) $Ad(k)(\am_i)=\am_i, \forall k\in K, i\in\{1,\cdots,r\}$;\\
(3) the map
$\varphi:K\times\am_1\times\cdots\times\am_r\rightarrow G$ defined
by $\varphi(k,X_1,\cdots,X_r)=ke^{X_1}\cdots e^{X_r}$ is a
diffeomorphism.\\
Now suppose $k_1,k_2\in K$ are conjugate in $G$, that is, there
exists some $g\in G$ such that $k_2=gk_1g^{-1}$. Rewrite this
equality as $k_1^{-1}gk_1=k_1^{-1}k_2g$, and write
$g=ke^{X_1}\cdots e^{X_r}$, where $k\in K$, $X_i\in\am_i$. Then we
have
$$
(k_1^{-1}kk_1)e^{Ad(k_1^{-1})X_1}\cdots
e^{Ad(k_1^{-1})X_r}=(k_1^{-1}k_2k)e^{X_1}\cdots e^{X_r},
$$
that is,
$$
\varphi(k_1^{-1}kk_1,Ad(k_1^{-1})X_1,\cdots,Ad(k_1^{-1})X_r)=
\varphi(k_1^{-1}k_2k,X_1,\cdots,X_r).
$$
Since $\varphi$ is a diffeomorphism, we have
$k_1^{-1}kk_1=k_1^{-1}k_2k$, that is, $k_2=kk_1k^{-1}$. Hence
$k_1,k_2$ are conjugate in $K$.
\end{proof}

\begin{corollary}\label{C:max}
Let $G$ be a connected Lie group, $K$ a maximal compact subgroup
of $G$. Then $i_*:\mathcal{E}_n(K)\rightarrow \mathcal{E}_n(G)$ is
a bijection.
\end{corollary}

\begin{proof}
Suppose $E^i\in \mathcal{E}_n(G)$, and choose an arbitrary element
$g\in E^i$. By Theorem 3.1 of Chapter XV in \cite{Ho}, there
exists $h\in G$ such that $hgh^{-1}\in K$. But $E^i$ is a
conjugacy class in $G$, so $hgh^{-1}\in E^i$. This proves $E^i\cap
K\neq\emptyset$. Hence $i_*$ is surjective.

The above argument in fact shows $E^i\cap K$ is a non-empty union
of some elements of $\mathcal{E}_n(K)$ for each $E^i$. We claim
that such a union consists of just one element of
$\mathcal{E}_n(K)$. In fact, if $E_K^1,E_K^2\in \mathcal{E}_n(K)$
such that $E_K^1\cup E_K^2\subset E^i\cap K$, choose $k_1\in
E_K^1, k_2\in E_K^2$. Since $k_1$ and $k_2$ belong to $K^i$, they
are conjugate in $G$. By Lemma \ref{L:max}, they are conjugate in
$K$. But  $E_K^1,E_K^2$ are conjugacy classes in $K$, this forces
$E_K^1=E_K^2$. In other word, the map $i_*$ is injective.
\end{proof}

The following conclusion is obvious from the proof the Corollary
\ref{C:max}.

\begin{corollary}
Let $G$ be a connected Lie group, $K$ a maximal compact subgroup
of $G$. Then the intersection of each element of
$\mathcal{E}_n(G)$ with $K$ is an elememt of
$\mathcal{E}_n(K)$.\qed
\end{corollary}

\begin{corollary}
Let $G$ be a connected Lie group. Then $\mathcal{E}_n(G)$ is a
finite set.
\end{corollary}

\begin{proof}
Let $K$ be a maximal compact subgroup of $G$, and let $T$ be a
maximal torus of $K$. By the theory of compact Lie groups, the
natural map $\mathcal{E}_n(T)\rightarrow \mathcal{E}_n(K)$ is
surjective. But $\mathcal{E}_n(T)$ is obvious finite, so
$\mathcal{E}_n(K)$ is also finite. By Corollary \ref{C:max},
$\mathcal{E}_n(G)$ and $\mathcal{E}_n(K)$ have the same
cardinalities. So $\mathcal{E}_n(G)$ is finite.
\end{proof}


\section*{Acknowledgments}

The first author would like to thank the following people:
Professor Z. Hajto and M. Shiota for valuable conversations,
Professor Jiu-Kang Yu for helpful discussion on algebraic groups
and for pointing out some errors in the first draft of the paper,
Professor K.-H. Neeb for informing the author a Lie-theoretic
proof of Theorem \ref{T:main} after the first draft of this paper
was available, as well as Dong Wang and Kuihua Yan for friendly
help.

\end{document}